
\baselineskip=14pt
\parskip=10pt

\font\eightrm=cmr8 

\magnification=\magstephalf
\def\A{{\cal A}}

\def\1{{\overline{1}}}
\def\2{{\overline{2}}}
\parindent=0pt
\overfullrule=0in

\def\frac#1#2{{#1 \over #2}}
\centerline
{\bf 
Using Symbolic Computation to Explore Generalized Dyck Paths and Their Areas
}
\bigskip
\centerline
{\it AJ BU and Doron ZEILBERGER}
\bigskip

\qquad {\it Dedicated to Symbolic Computation pioneer Bruno Buchberger (b. 22 October, 1942) on his $80^{th}$ birthday}

{\bf Abstract}: We show the power of Bruno Buchberger's seminal Gr\"obner Basis algorithm, interfaced, seamlessly, with what we call
{\it symbolic dynamical programming}, to {\it automatically} generate algebraic equations satisfied by the generating functions
enumerating so-called {\it Generalized Dyck Walks}, i.e. 2D walks that start and end on the $x$-axis, and never dip below it, for an {\it arbitrary}
set of steps. More impressively, we combine it with  calculus (that Maple knows very well!), to
automatically compute generating functions for the sum-of-the-areas of these generalized Dyck paths,
and even for the sum of any given power of the areas, enabling us to get statistical information about the area under a
random generalized Dyck path.

{\bf Maple package and Output Files}

This article is accompanied by a Maple package, {\tt GDW.txt}, and numerous output files, some of which will be referred to later.

They are all obtainable from the front of this article

{\tt https://sites.math.rutgers.edu/\~{}zeilberg/mamarim/mamarimhtml/area.html } \quad .

{\bf Bruno Buchberger: the Gauss of Non-Linear Algebra}

Using {\it Gaussian elimination}, we can solve, efficiently, any system of {\bf linear equations} with any (finite) number of unknowns.
Using the {\bf Buchberger algorithm} [B1] [B2] (see the modern classic [CLO] for a lucid and engaging account) we can solve
any system of {\bf non-linear equations} with any (finite) number of unknowns. Alas, due to the non-linearity, we can only
go so far, but yet, with modern computers, and {\it computer algebra systems} (our favorite being Maple), one can do a lot.

In the applications that we need, the scenario is that we have $N$ quantities $X_1, \dots,  X_N$ and
they satisfy $N$ polynomial equations
$$
P_i(t;X_1, \dots, X_n)=0 \quad, \quad 1 \leq i \leq N \quad,
$$
where $t$ is variable (or multi-variable), that may be viewed as a parameter. We are {\bf really} only interested in the first quantity, $X_1=X_1(t)$.

In our application, these equations are derived using combinatorial considerations, that we can teach to the computer.
There is one {\bf primary} object, $X_1$,  that we really care about, but in order to determine it, we need to introduce {\bf secondary} objects,
$X_2, \dots, X_n$, that we don't care about, and whose sole {\it raison d'\^etre} is that they enable us to get a grip on $X_1$.

All we want is a {\bf succinct} algebraic relation of the form $Q(X_1,t)=0$, where
$Q$ is a bi-variate polynomial of $X_1$ and $t$.

In order to {\bf eliminate} the unwanted quantities $X_2, \dots , X_n$,  and
get a pure equation only involving $X_1$ (and the auxiliary variable (or variables), $t$),  the command in {\tt Maple} is:

$$
Groebner[Basis](\{P_1,P_2, \dots, P_n \} ,plex(X_2,X_3, \dots, X_n, X_1))[1]; \, .
$$

To take a random example (we just made up), suppose that we have the system of three equations in the three unknowns $X,Y,Z$ that
depend on some auxiliary variable $t$:
$$
X=Y^2+Z^2+1 \quad, \quad Y=X^2+3Z^2+t \quad, \quad Z=XYZ+t+1 \quad,
$$
and we want to {\bf eliminate} $Y$ and $Z$, i.e. get an algebraic equation in $X=X(t)$ that only involves, in addition to $X$, the variable $t$. Then type in a Maple session

{\tt Groebner[Basis]( $\{$ X-Y**2-Z**2-1,Y-X**2-3*Z**2-t, Z-X*Y*Z-t-1 $\}$ ,plex(Y,Z,X))[1]=0;} 

getting that our quantity of interest $X=X(t)$, satisfies the following degree-$12$  polynomial equation in $X(t)$ (abbreviated $X$), with coefficients that are polynomials of $t$:
$$
X^{12}+6 X^{11}+\left(4 t +3\right) X^{10}+\left(18 t -21\right) X^{9}+\left(6 t^{2}-8\right) X^{8}+\left(18 t^{2}-44 t -9\right) X^{7}
$$
$$
+\left(4 t^{3}+9 t^{2}+20 t +60\right) X^{6}+\left(6 t^{3}+29 t^{2}+96 t +43\right) X^{5}+\left(t^{4}+30 t^{3}+16 t^{2}-14 t -24\right) X^{4}
$$
$$
+\left(52 t^{3}+177 t^{2}+232 t +50\right) X^{3}
+\left(18 t^{4}-12 t^{3}+89 t^{2}+242 t +149\right) X^{2}
$$
$$
+\left(54 t^{3}+12 t^{2}-150 t -105\right) X +81 t^{4}+378 t^{3}+612 t^{2}+396 t +99 \, = \, 0 \quad.
$$

{\bf The Art of Counting in General}

In {\it enumerative combinatorics}, we  are interested in {\it counting}  families of sets $\{A(n)\}$, 
defined by some combinatorial conditions, indexed by (one or more) discrete parameter(s) $n$.
The best possible scenario is an {\bf explicit} (aka {\it closed form}) expression for 
$a(n):=|A(n)|$, where as usual $|S|$ denotes the number of elements of the finite set $S$.

For example, the number of words in the alphabet $\{0,1\}$ of length $n$ (Ans.: $2^n$), or the number of permutations of length $n$ (Ans.: $n!$).

Alas, in many cases there is no {\it closed form} (unless you consider $\sum_{s \in A(n)} 1$ as one), and one looks for a closed-form, or failing this, an {\it equation}, satisfied by
the (ordinary) generating function
$$
X(t) \, = \, \sum_{n=0}^{\infty} a(n)t^n \quad.
$$

In many applications, including all those in the present paper, it is better to think of $X(t)$ as a {\it weight enumerator} of the (usually) `infinite' set
$$
\A \, :=\, \bigcup_{n=0}^{\infty} A(n) \quad, 
$$
and define, $Weight(s):=t^n$ if $s\in A(n)$. The generating function, $X(t)$, is the {\bf total weight} (i.e. sum of all the weights of its members) of $\A$, denoted by $Weight(\A)$.

{\bf Reminders about Counting Classical Dyck and Motzkin Walks}

In order to illustrate our general approach, already initiated in [AyZ], and further pursued in [TZ], [EK1], and [EK2], let's revisit the very classical problem
of counting (ungeneralized) {\bf Dyck walks}. These are walks that start at the origin $(0,0)$, end somewhere on the $x$-axis, 
using the {\bf atomic steps} $[1,1],[1,-1]$ (i.e. moving from $(x,y)$ to either $(x+1,y+1)$ or $(x+1,y-1)$, respectively)
and {\bf crucially}, always stay in the upper-half plane $ y \geq 0$. 

The weight of a walk $w$ is $t^{NumberOfSteps(w)}$, in other words $t^n$, if the endpoint of our walk is $(n,0)$.

Let $f[0,0](t)$ be our object of desire, the weight-enumerator of all Dyck walks.  There are two possibilities.

$\bullet$ The walk  is empty (of length $0$), with weight $t^0=1$.

$\bullet$ The walk, let's call it $w$, sooner of later, again touches the $x$-axis. Let $w_1$ be the prefix of $w$ consisting of this portion and let $w_2$ be the rest, so we can
express $w$ as the {\bf catenation} $w=w_1\,w_2$, where $w_2$ is a (possibly empty) shorter walk weight-counted by $f[0,0](t)$, but $w_1$ is a {\bf strict Dyck path},
a walk that starts and ends on the $x$-axis, and except at the endpoints, is {\bf strictly} above the $x$-axis. Let's call the weight-enumerator of
these kind of walks $g[0,0](t)$. 

Taking weights, we get a first equation
$$
f[0,0](t)=1 + g[0,0](t) \cdot f[0,0](t) \quad .
$$

Note that we were {\bf forced} to introduce the quantity $g[0,0](t)$.

Let's look at the anatomy of such a walk that is weight-enumerated by $g[0,0](t)$. Since it is non-empty, 
and can't get below the $x$-axis, it {\bf must} start with the step $[1,1]$, i.e. the first step is from the
origin $(0,0)$ to $(1,1)$. Since sooner or later it must return to the $x$-axis, its last step must be
a down step $[1,-1]$, from $(n-1,1)$ to $(n,0)$ (say). Removing this first and last step is
a walk that is shorter by two steps that is {\it weakly} above the line $y=1$. But these
walks are in obvious bijection with walks weight-counted by $f[0,0](t)$. Hence we get a {\bf second} equation
$$
g[0,0](t)\,=\, t^2 \cdot f[0,0](t)\ \quad .
$$
So we get a system of two equations and two unknowns.

In this case, we don't need a computer, or Bruno Buchberger, to solve this system, and eliminating $g[0,0](t)$ gives the
quadratic equation for $f[0,0](t)$ (let's abbreviate it to $f[0,0]$)
$$
t^2\,f[0,0]^2- f[0,0]+1 \, = \, 0 \quad ,
$$
that thanks to the Babylonians can be solved in terms of radicals, yielding the Catalan numbers.

Similarly, for {\bf Motzkin walks}, where the set of steps is $\{ [1,1],[1,0],[1,-1] \}$, we still only have one extra quantity $g[0,0](t)$, and the two equations are
$$
\{f[0,0]=1 \,+ \, t f[0,0]\, + \, g[0,0]\cdot f[0,0]  \quad, \quad g[0,0]\,=\, t^2 \cdot f[0,0]\} \quad, 
$$
that again lead to a simple quadratic equation, whose solution is the generating function for Motzkin numbers.

The primary novelty of the present article is to keep track of the {\bf area} statistics. Let $f[0,0](t,q)$, that we will still call $f[0,0](t)$ for short
(not to be confused with the previous $f[0,0](t)$), but it is implied that there is also a variable $q$, that may be viewed as a parameter,
be the weight-enumerator of Dyck walks according to the bi-variate weight
$$
Weight(w):=t^{LengthOf(w)}\, q^{AreaUnder(w)} \quad.
$$
One way to define the `area' without geometry is as the sum of the $y$-coordinates of all the intermediate points on the walk.

Going back to the classical Dyck  case, instead of a system of {\it algebraic equations}, we get a system of
{\it functional equations}. It is easy to see that we have
\quad .
$$
\{f[0,0](t) \,= \, 1 + g[0,0](t) \cdot f[0,0](t) \quad ,   \quad g[0,0](t) \,= \, qt^2 \cdot f[0,0] (qt) \} \quad .
$$
Of course one can eliminate $g[0,0](t)$ and get a pure functional equation for $f[0,0](t)$ (let's call it $X(t)$).
$$
X(t)=1 \,+ \,qt^2 \cdot X(t) \cdot X(qt) \quad.
$$
Alas $f[0,0](t)=X(t)$ is no longer an algebraic formal power series.

The study of the area under Dyck paths was pioneered in [MSV] and further explored in
[Wo],[SRW], and for higher moments by Robin Chapman [C]. We will soon see how to fully automate this  to a much more general
setting. In particular, our Maple package, {\tt GDW.txt}, can confirm, in a few nano-seconds, all their beautiful human-generated results.

Similarly, for the Motzkin case, the system of functional equations is:
$$
\{
f[0,0](t)=1 \,+\, t \cdot f[0,0](t) \,+ \, g[0,0](t) \cdot f[0,0](t) \quad ,   \quad g[0,0](t)\,=\,qt^2 \cdot f[0,0] (qt) \} \quad ,
$$
and the pure functional equation for $f[0,0](t)$ (again let's call it $X(t)$) is:
$$
X(t)=1 \,+ \, t\cdot X(t) \,+ \,qt^2\cdot X(t) \cdot X(qt) \quad.
$$

The bi-variate generating function contains lots of statistical information. To get the average area under a random Dyck walk, we must first compute
the generating function for the  `sum of the areas of all walks' that is given by 
$$
\frac{d}{dq}f[0,0](t) {\big \vert}_{q=1} \quad,
$$

(recall that $f[0,0](t)$ is really a function of both $t$ and $q$.). More generally, the generating function for the sum of the $r$-th powers of the areas, is
$$
\left ( q\frac{d}{dq} \right)^r f[0,0](t) {\big \vert}_{q=1} \quad.
$$

If you are interested in these quantities up to the $r$-th power of the area, we can taylor expand $f[0,0](t)$ (and similarly $g[0,0](t)$)
as a Taylor expansion about $q=1$:
$$
f[0,0](t)= f[0,0,0](t)+f[0,0,1](t)(q-1)+f[0,0,2](t)(q-1)^2+ \dots + f[0,0,r](t) (q-1)^r+ \dots
$$
where $f[0,0,0](t),f[0,0,1](t),\dots$ are formal power series of $t$ alone, and do not depend on $q$. 
(Note that what we now call $f[0,0,0](t)$ is identical to what we called $f[0,0](t)$ in the original, straight-enumeration case).

It would follow from the algorithm
to be described shortly, that these are all {\bf algebraic} formal power series, and thanks to the {\bf Buchberger Gr\"obner basis algorithm}, we can actually find the
{\it pure} equations that they each satisfy.

Plugging this {\it series expansion} into the functional equation, would introduce quantities like $f[0,0,0](qt)$.

For this we need an elementary lemma from Calculus:

{\bf Simple Lemma:} Let $f(t)$ be a formal power series of a single variable $t$, and $q$ be another variable, then 
$$
f(qt) \, = \, f(t) + (q-1)tf'(t)+ \frac{1}{2} (q-1)^2 t^2f''(t) + \dots + \frac{1}{r!} (q-1)^r t^r f^{(r)}(t) + \dots \quad .
$$

So we can Taylor-expand it with respect to $q$ around $q=1$, and thanks to Calculus, that Maple knows so well, {\it automatically}
express it in terms of $f[0,0,0]'(t), f[0,0,1]'(t) \dots \quad$. It seems that we have too many quantities, and that we need more
equations. But we can also (automatically!) differentiate with respect to $t$ the algebraic equation satisfied by $f[0,0,0](t)$, 
(and if needed, by $f[0,0,1](t)$ etc.) using
{\it implicit differentiation} (that Maple also knows!),  and collect terms, and at the {\bf end of the day} we
would have as many algebraic equations as quantities. If we want to focus, say, on $f[0,0,1](t)$ (alias $(\frac{\partial}{\partial q})f[0,0](t){\big \vert}_{q=1}$,
alias the generating function for the `sum of areas' of the walks), we can use Buchberger's algorithm as above with {\it plex}, singling out $f[0,0,1](t)$. Ditto for higher (factorial) moments, i.e. $f[0,0,2](t), f[0,0,3](t), \dots$
except that things do get more and more complicated, even for computerkind.

{\bf Straight Enumeration of Generalized Dyck Paths}

Let's recall the method initiated in [AyZ], extended in [EK1] and [EK2], and continued in [TZ], but with a new implementation, that makes full use of Gr\"obner bases.
Alternative approaches to the problem of (straight) enumeration have been undertaken by Banderier et. al. [BKKKKNW], using the {\it kernel method}, and
in [EkhZ], using {\it ``Guess and Check''}.

The {\bf input} is an arbitrary set of integers, $S$, and our goal is to find a pure algebraic equation of the form
$$
\sum_{i=0}^d p_i(t) X(t)^i \quad=0,
$$
where $p_i(t)$ are polynomials in $t$ for the following formal power series
$$
X(t)= \sum_{n=0}^{\infty} a_S(n)\, t^n  \quad,
$$
where $a_S(n)$ is the number of walks from $(0,0)$ to $(n,0)$, in the 2D lattice, with the set of allowable steps
$$
(x,y) \rightarrow (x+1, y+s) \quad, \quad s \in S \quad,
$$
and that always stay in the upper half plane $y \geq 0$. Equivalently, and more useful for us, it is the {\it weight enumerator}, according to the weight
$t^{NumberOfSteps}$, of the (infinite) set of all these walks.

{\bf Remark}: If you don't insist on always staying above the $x$-axis, then the number of such $n$-step walks from the origin to the $x$-axis, is the constant term of
$$
\left ( \sum_{s \in S} t^s \right )^n \quad \quad,
$$
and it is possible, very fast, to use the {\bf Almkvist-Zeilberger algorithm} [AlZ]  (see [D] for a great exposition)
to get a linear recurrence equation with polynomial coefficients.

Note that if $S$ only consists of positive integers, or only consists of negative integers, then $X(t)$ is trivially $1$. So for things to be non-trivial,
$S$ must have at least one positive member and at least one negative member.

We will rename $X(t)$, $f[0,0](t)$, since we need to introduce auxiliary quantities $f[a,b](t)$ and $g[a,b](t)$.

For integer $a \geq 0$ and $b \geq 0$,

$\bullet$ Let $f[a,b](t)$ be the weight-enumerator of walks with a set of steps given by $S$,  that start at the point $(0,a)$ and end on the horizontal line $y=b$
and stay {\it weakly} above the $x$-axis.

$\bullet$ Let $g[a,b](t)$ be the weight-enumerator of  non-empty walks with a set of steps given by $S$,  that start at the point $(0,a)$ and end on the horizontal line $y=b$
and stay {\it strictly} above the $x$-axis, except at an endpoint when $a=0$ or $b=0$. If both $a>0$ and $b>0$ then we have no need for it, and we declare that it is $0$.

For each $f[a,b](t)$ and $g[a,b](t)$ that will show up we would need to set its own equation.

{\bf How to form the Equations for the $f[a,b]$?}

$\bullet$ If $a>0$ and $b>0$ then

$$
f[a, b] \,= \,g[a, 0] \cdot f[0, b] \,+ \,f[a - 1, b - 1] \quad .
$$

{\bf Explanation:} if such a walk, that starts at $y=a$ and ends at $y=b$, touches the $x$-axis, then the first portion until that first encounter is a walk weight-counted by $g[a,0]$, while the
second part is a walk weight-counted by $f[0,b]$ (it starts at $y=0$, ends at $y=b$ and stays weakly above the $x$-axis). On the other hand, if it never touches
the $x$-axis it must stay in $y \geq 1$, so lowering it by $1$ unit, yields a walk weight-counted by $f[a-1,b-1]$.

$\bullet$ If $a>0$ and $b=0$ then
$$
f[a, 0] \,= \, g[a,0] \cdot f[0, 0] \quad .
$$
{\bf Explanation:} A walk that starts at $(0,a)$ and ends on $y=0$, must meet the $x$-axis for the first time, this initial portion is weight-counted by $g[a,0]$, the remaining portion of the walk is weight-counted by $f[0,0]$.

Similarly

$\bullet$ If $a=0$ and $b>0$ then
$$
f[0, b] \,= \, f[0,0] \cdot g[0, b] \quad .
$$

$\bullet$ If $a=0$ and $b=0$ and $0 \not \in S$
$$
f[0, 0] \,= \, 1 \,+ \, g[0,0] \cdot f[0, 0] \quad .
$$
{\bf Explanation:} A walk that starts and ends on $y=0$ may be the empty walk (weight=$1$).
If not, it can be broken up into the portion until the first encounter with the $x$ axis, that is weight-counted by $g[0,0]$, followed by any-old (shorter) walk weight-counted by $f[0,0]$.

$\bullet$ If $a=0$ and $b=0$ and $0 \in S$
$$
f[0, 0] \,= \,  1 \,+\, t\cdot f[0,0] \,+ \, g[0,0] \cdot f[0, 0] \quad .
$$
{\bf Explanation:} Now the first step could also be $(0,0) \rightarrow (1,0)$, and after that it is a typical walk weight-counted by $f[0,0]$.

{\eightrm This is implemented in procedure {\tt MakeEqF(f,g,x,a,b,S)} in our Maple package. In fact, it is more general, keeping track of the individual steps.}

{\bf How to form the Equations for the $g[a,b]$?}

We also need to set up equations for $g[a,b]$ for those $(a,b)$ that would be needed.

Let $P$ be the subset of $S$ consisting of the (strictly) positive members of $S$, and
let $N$ be the subset of $S$ consisting of the (strictly) negative members of $S$, so if $0 \in S$ then
$$
S= P \cup N \cup \{0\} \quad,
$$
while,  if $0 \not \in S$ then
$$
S= P \cup N  \quad .
$$

$\bullet$ If $a=0$ and $b>0$ then
$$
g[0, b] =  t \left ( \sum_{i \in P} f[a+i-1,b-1] \right ) \quad .
$$
{\bf Explanation:} Such a walk must start with a step of the form $[1,i]$ where $i \in P$, and then it is always in $y \geq 1$. Removing that first step (weight $t$) and `lowering' it by $1$ unit, is in
bijection with a walk weight-counted by $f[a+i-1,b-1]$.

$\bullet$ If $a>0$ and $b=0$ then
$$
g[a, 0] =  t \left ( \sum_{j \in N} f[a-1,b-j-1] \right ) \quad .
$$
{\bf Explanation:} Such a walk must end with a step of the form $[1,j]$ where $j \in N$, and before that it is always in $y \geq 1$. Removing that last step (weight $t$) and `lowering' it by $1$ unit, is in
bijection with a walk weight-counted by $f[a-1,b-j-1]$.

$\bullet$ If $a=0$ and $b=0$ then
$$
g[0,0]= t^2 \left ( \sum_{i \in P} \sum_{j \in N} f[a+i-1,b-j-1] \right ) \quad.
$$
{\bf Explanation:} Every non-empty walk that starts at $y=0$ and ends at $y=0$ and that never touches the $x$-axis except at the two endpoints, must start with a positive step $[1,i]$ and
end with a negative step $[1,j]$. Deleting the first and last steps gives you a shorter walk (with two steps shorter) that is always in the region $y \geq 1$.
Lowering the remaining walk  by $1$ unit gives a walk weight-counted by $f[a+i-1,b-j-1]$.

{\eightrm This is implemented in procedure {\tt MakeEqG(f,g,x,a,b,S)} in our Maple package. Again, it is more general keeping track of the individual steps.

To get the full system of equations with a set of steps {\tt S}, type:

MakeSysT(f,g,t,S);

It returns the set of equations, followed by the set of quantities that participate. For example

{\tt MakeSysT(f,g,t, $\{$ 1,2,-1,-2 $\}$ )[1];} gives

$$
\{f_{00} = f_{00} g_{00}+1 \quad, \quad f_{01} = f_{00} g_{01}
\quad, \quad f_{10} = g_{10} f_{00} \quad , \quad f_{11} = f_{01} g_{10}+f_{00}, 
$$
$$
g_{00} = 
t^{2} f_{00}+t^{2} f_{01}+t^{2} f_{10}+t^{2} f_{11} \quad , \quad g_{01} = t f_{00}+t f_{10} \quad , \quad  g_{10} = t f_{00}+t f_{01}\} \quad,
$$

while {\tt MakeSysT(f,g,t, $\{$1,2,-1,-2 $\}$)[2];} gives the set of quantities
$$
\{f_{00} \quad , \quad f_{01} \quad , \quad f_{10} \quad , \quad f_{11} \quad, \quad g_{00} \quad , \quad g_{01} \quad, \quad g_{10}\} \quad .
$$
}

{\bf  Weighted Enumeration of Generalized Dyck Paths According to the Area}

We can teach Maple how, {\it all by itself}, set up a system of {\bf functional equations}, for the $q$-analog, where one keeps track of the area.
Everything is analogous.

For integer $a \geq 0$ and $b \geq 0$,

$\bullet$ Let $f[a,b](t)$ be the weight-enumerator of  walks with a set of steps given by $S$,  that start at the point $(0,a)$ and end on the horizontal line $y=b$
and stay {\it weakly} above the $x$-axis, but now the {\it weight} of a walk is not just $t^{length}$ but rather $t^{length}\cdot q^{AreaUnder}$.
For the sake of notational convenience, We still write $f[a,b](t)$ rather than $f[a,b](t,q)$, but one should keep in mind that they all
depend on $q$. We can think of $q$ as a parameter.

Similarly,

$\bullet$ Let $g[a,b](t)$ be the weight-enumerator (in the above sense of $t^{length}\cdot q^{AreaUnder}$)
of non-empty walks with a set of steps given by $S$,  that start at the point $(0,a)$ and end on the horizontal line $y=b$
and stay {\it strictly} above the $x$-axis, except at an endpoint when $a=0$ or $b=0$. If both $a>0$ and $b>0$ then we have no need for it and we declare that it is $0$.

{\bf How to form the Functional Equations for $f[a,b](t)$?}

$\bullet$ If $a>0$ and $b>0$ then
$$
f[a, b](t) = g[a, 0](t) \cdot f[0, b](t) \,+ \,f[a - 1, b - 1](q\,t) \quad .
$$

{\bf Explanation:} if such a walk that starts at $y=a$ and ends at $y=b$ touches the $x$-axis, then the first portion until that first encounter is a walk weight-counted by $g[a,0](t)$, while the
second part is a walk weight-counted by $f[0,b](t)$ (it starts at $y=0$, ends at $y=b$ and stays weakly above the $x$-axis). On the other hand, if it never touches
the $x$-axis, it must stay in $y \geq 1$, so lowering it by $1$ unit, results in a walk weight-counted by $f[a-1,b-1](t)$. {\bf But}, by lowering it, we lost some area!.
If the length of the walk is $n$, then we lost $n$ units of area (the area of an $n \times 1$ rectangle), hence this is bi-weight-enumerated  by $f[a-1,b-1](q\,t)$,
(rather than $f[a-1,b-1](t)$.)

$\bullet$ If $a>0$ and $b=0$ then
$$
f[a, 0](t) \,= \, g[a,0](t) \cdot f[0, 0](t) \quad .
$$
{\bf Explanation:} Each walk that starts at $y=a$ must meet the $x$-axis for the first time, this portion is weight-counted by $g[a,0](t)$, the remaining portion of the walk is weight-counted by $f[0,0](t)$.

Similarly:

$\bullet$ If $a=0$ and $b>0$ then
$$
f[0, b](t) \,= \, f[0,0](t) \cdot g[0, b](t) \quad,
$$

$\bullet$ If $a=0$ and $b=0$ and $0 \not \in S$
$$
f[0, 0](t) \,= \, 1 \, + \, g[0,0](t) \cdot f[0, 0](t) \quad .
$$
{\bf Explanation:} Each walk that starts and ends on $y=0$ must either be the empty walk (of length $0$, area $0$, and hence with weight $t^0 \cdot q^0=1$), or else
meet the $x$-axis for the first time (weight-counted by $g[0,0](t)$), followed by any-old walk (possibly empty) weight-counted by $f[0,0](t)$.

$\bullet$ If $a=0$ and $b=0$ and $0 \in S$
$$
f[0, 0](t) \,= \, 1 \,+ \, t \cdot f[0,0](t) \,+ \, g[0,0](t) \cdot f[0, 0](t) \quad .
$$
{\bf Explanation:} Now the first step could also be $(0,0) \rightarrow (1,0)$, and after that it is a typical walk weight-counted by $f[0,0](t)$, and there is no area gain.

{\eightrm This is implemented in procedure {\tt qMakeEqF(f,g,t,q,a,b,S)} in our Maple package. }

{\bf How to form the Functional Equations for $g[a,b](t)$?}

We also need to set up equations for $g[a,b](t)$ for those $(a,b)$ that would be required.

As before:

Let $P$ be the subset of $S$ consisting of the (strictly) positive members of $S$, and
let $N$ be the subset of $S$ consisting of the (strictly) negative members of $S$, so if $0 \in S$ then
$$
S= P \cup N \cup \{0\} \quad,
$$
while,  if $0 \not \in S$ then
$$
S= P \cup N  \quad .
$$

$\bullet$ If $a=0$ and $b>0$ then
$$
g[0, b](t) =  t \left ( \sum_{i \in P} q^{i/2}\,f[a+i-1,b-1](qt) \right ) \quad .
$$
{\bf Explanation:} Such a walk must start with a step of the form $[1,i]$ where $i \in P$, and then it is always in $y \geq 1$. Removing that first step reduces the area by $1\times i/2 =i/2$
and `lowering' by $1$ unit, it looses $n$ units of area (if the remaining path has length $n$), to account for the lost area we need $f[a+i-1,b-1](qt)$, (rather than  $f[a+i-1,b-1](t)$).

Similarly:

$\bullet$ If $a>0$ and $b=0$ then
$$
g[a, 0](t) =  t \left ( \sum_{j \in N} q^{-j/2} f[a-1,b-j-1](qt) \right ) \quad .
$$

$\bullet$ If $a=0$ and $b=0$ then
$$
g[0,0](t)= t^2 \left ( \sum_{i \in P} \sum_{j \in N} q^{i/2-j/2}\,f[a+i-1,b-j-1](qt) \right ) \quad.
$$
{\bf Explanation:} Every walk that starts at $y=0$ and ends at $y=0$ and that never touches the $x$-axis except at the endpoints, must start with a positive step $[1,i]$ and
end with a negative step $[1,j]$. Deleting the first and last steps gives you a shorter walk (with two steps shorter) that is always in the region $y \geq 1$.
Removing the first step (i.e. $(0,0) \rightarrow (1,i)$)  reduces the area by $i/2$.
Removing the last step (i.e. $(n-1,-j) \rightarrow (n,0)$)  reduces the area by $-j/2$.

Lowering it by $1$ unit gives a walk weight-counted by $f[a+i-1,b-j-1](qt)$.

{\eightrm This is implemented in procedure {\tt qMakeEqGt(f,g,t,q,a,b,S)} in our Maple package.}

To get the full system of functional equations, followed by the quantities that feature in them,  with set of steps {\tt S}, type:

{\tt qMakeSysT(f,g,t,q,S);} \quad.

For example, if the set of steps is $S=\{2,1,0,-1,-2\}$, typing

{\tt qMakeSyst(f,g,t,q, $\{$ 2,1,0,-1,-2 $\}$ )[1];}

gives

$$
\{
f_{01}\! \left(t \right)-f_{00}\! \left(t \right) g_{01}\! \left(t \right) \quad ,  \quad
f_{10}\! \left(t \right)-g_{10}\! \left(t \right) f_{00}\! \left(t \right) \quad, 
$$
$$
f_{11}\! \left(t \right)-f_{01}\! \left(t \right) g_{10}\! \left(t \right)-f_{00}\! \left(q t \right) \quad ,  \quad
g_{01}\! \left(t \right)-t \sqrt{q}\, f_{00}\! \left(q t \right)-t q f_{10}\! \left(q t \right) \quad, 
$$
$$
g_{10}\! \left(t \right)-t q f_{01}\! \left(q t \right)-t \sqrt{q}\, f_{00}\! \left(q t \right),  \quad
f_{00}\! \left(t \right)-f_{00}\! \left(t \right) g_{00}\! \left(t \right)-f_{00}\! \left(t \right) t -1 \quad,
$$
$$
g_{00}\! \left(t \right)-t^{2} q^{\frac{3}{2}} f_{01}\! \left(q t \right)-t^{2} q f_{00}\! \left(q t \right)-t^{2} q^{2} f_{11}\! \left(q t \right)-t^{2} q^{\frac{3}{2}} f_{10}\! \left(q t \right)
\} \quad ,
$$
while to see the set of featured quantities type

{\tt qMakeSyst(f,g,t,q, $\{$ 2,1,0,-1,-2 $\}$)[2];} , getting
$$
\{f_{00}\! \left(t \right) \quad , \quad f_{00}\! \left(q t \right) \quad , \quad
f_{01}\! \left(t \right) \quad , \quad f_{01}\! \left(q t \right) \quad ,  \quad
f_{10}\! \left(t \right) \quad , \quad f_{10}\! \left(q t \right) \quad ,  \quad
$$
$$
f_{11}\! \left(t \right) \quad , \quad f_{11}\! \left(q t \right) \quad ,  \quad
g_{00}\! \left(t \right) \quad , \quad g_{01}\! \left(t \right) \quad ,  \quad
g_{10}\! \left(t \right)\} \quad .
$$

{\bf From Functional Equations to Algebraic Equations}

After the computer finds the system of {\it functional} equations described above, we instruct it to find a system {\it algebraic} equations
for the `components' of the $f[a,b](t)$ (and we also need $g[a,b](t)$). If we are interested in finding the pure algebraic equation satisfied by
the formal power series $f[0,0,0](t), f[0,0,1](t), \dots, f[0,0,k](t)$. We write
$$
f[0,0](t)=\sum_{i=0}^{k} f[0,0,i](t) (q-i)^i + O((q-1)^{k+1}) \quad,
$$
and analogously for the other $f[a,b](t)$ and $g[a,b](t)$, then expand in powers of $(q-1)$, collect terms, use the Lemma, and get more
equations by differentiating with respect to $t$ each of these equations up to the $k$-th derivative, using implicit differentiation.

Because of the extreme complexity, we decided only to implement this scheme for $k=1$, i.e. for finding the
algebraic equation satisfied by the generating function for the `sum of the areas'.

{\eightrm

This is implemented in procedure
{\tt qEqGFt(S,X,t)} . For example, to get the algebraic equation for the generating function for `sum of areas' of the classical Dyck paths, type:

{\tt qEqGFt($\{$1,-1$\}$,X,t); }

getting
$$
t^2 \,- \, (4t^2-1)(2t^2-1)X \,+ \, t^2(4t^2-1)^2\,X^2 \,=\,0 \quad .
$$

(This is A8549 of [Sl], {\tt https://oeis.org/A008549}).

For Motzkin walks, typing

{\tt qEqGFt($\{$1,0,-1$\}$,X,t);},  gives
$$
t^{2}-\left(3 t -1\right) \left(t +1\right) \left(t^{2}+2 t -1\right) X +t^{2} \left(3 t -1\right)^{2} \left(t +1\right)^{2} X^{2} \, = \, 0 \quad .
$$

(This is A57585 of [Sl], {\tt https://oeis.org/A057585}).

For a more complicated example, to get the pure algebraic equation satisfied by
the generating function for the `sum of the areas under generalized Dyck paths with set of steps $\{[1,2],[1,1],[1,0],[1,-1],[1,-2] \}$', type:

{\tt qEqGFt($\{$ 2,1,0,-1,-2 $\}$,X,t);}, getting, after less than a minute,
$$
t^{2} \left(775 t^{4}-1460 t^{3}+1006 t^{2}-264 t +24\right)
$$
$$
+\left(t -1\right) \left(5 t -1\right) \left(425 t^{6}-1520 t^{5}+1527 t^{4}-68 t^{3}-282 t^{2}+88 t -8\right) X 
$$
$$
-t \left(150 t^{5}+540 t^{4}-889 t^{3}-240 t^{2}+228 t -32\right) \left(t -1\right)^{2} \left(5 t -1\right)^{2} X^{2}
$$
$$
-2 t^{2} \left(5 t +4\right) \left(5 t^{3}-t^{2}-17 t +4\right) \left(t -1\right)^{3} \left(5 t -1\right)^{3} X^{3}
+t^{4} \left(5 t +4\right)^{2} \left(t -1\right)^{4} \left(5 t -1\right)^{4} X^{4}
\, = \, 0 \quad.
$$
}

(This is not (yet, May 15, 2023) in the OEIS. Note that the straight enumeration version is A104184 of [Sl], {\tt https://oeis.org/A104184}). For the sake of the OEIS, here are the first $30$ terms:
$$
0, 0, 3, 18, 113, 636, 3487, 18656, 98429, 514012, 2664690, 13737758, 70522801, 360806214, 1840913908, 
$$
$$
9371761174, 47621259557, 241601881822, 1224111502194, 6195045902854, 
$$
$$
31321134873744, 158217553824544, 798622703316154, 4028438371631942, 
$$
$$
20308239308212037, 102323623873153810, 515313296262175206,  2594054240062008690, 
$$
$$
13053194513626873348, 65659889953142043376 \quad .
$$

{\bf Strict Generalized Dyck paths}

If we want to count strict generalized Dyck paths (respectively, sum of the areas), i.e. paths that never touch the $x$-axis except at the endpoints, use
procedures {\tt EqGFtS(S,X,t)} and {\tt qEqGFtS(S,X,t)} respectively.

For the algebraic equation for the generating function for the sum of the areas under {\bf strict} classical Dyck paths, type

{\tt qEqGFtS($\{$1,-1$\}$,X,t);}, getting that the equation is
$$
\left(4 t^{2}-1\right) X +t^{2} \, = \, 0  \quad,
$$
that implies
$$
X(t)=\frac{t^{2}}{1-4 t^{2}} \quad,
$$
confirming, {\it purely automatically}, the following elegant proposition first discovered, and proved, in [SRW] (see also [C]):

{\bf Proposition} (Shapiro, Rogers, and Woan) The sum of the areas of the strict Dyck paths of length $2n$ is $4^{n-1}$.

What about sum-of-the-areas of strict Motzkin paths? Typing

{\tt qEqGFtS($\{$ 1 ,0,-1 $\}$,X,t);} gives
$$
\left(3 t^{2}+2 t -1\right) X +t^{2} \, = \, 0 \quad, 
$$
implying that
$$
X(t)=\frac{t^{2}}{1-2 t-3t^2} \quad.
$$
This is A015518[n-1] of [Sl] (see {\tt https://oeis.org/A015518}). This sequence has numerous combinatorial interpretations, but so far, the
connection to the sum of the areas under strict Motzkin paths escaped notice.

{\bf Efficient Computation of many terms}

It is well-known and fairly easy to see (e.g. the modern classic [KP]) that every algebraic formal power series is $D$-finite, and equivalently,
the sequence itself satisfies some {\bf linear recurrence equation with polynomial coefficients}. The Maple package {\tt gfun} 
(designed by Bruno Salvy and Paul Zimmermann) can do it for you, but it is just as easy to use {\it guessing}. Using {\it dynamical programming}
one can easily crank out many terms, and then use `guessing' (implemented in {\tt listtorec} in {\tt gfun}, but we have our own
home-made version). Now we can even talk about sum-of-square-of-areas, and even higher powers.

{\bf Sample Output}

{\bf Straight Enumeration of Generalized Dyck Paths}

$\bullet$ 
If you want to see $16$ theorems stating algebraic equations satisfied by generating functions enumerating generalized Dyck paths with set of steps that are subsets of 
$\{-2,-1,0,1,2\}$ (excluding trivial cases) and sometimes linear recurrences with polynomial coefficients satisfied by the sequences themselves, enjoy

{\tt https://sites.math.rutgers.edu/\~{}zeilberg/tokhniot/oGDW1.txt } \quad .

$\bullet$ 
If you want to see $94$ theorems stating algebraic equations satisfied by generating functions enumerating generalized Dyck paths with set of steps that are subsets of 
$\{-3,-2,-1,0,1,2,3\}$ (excluding trivial cases) and sometimes linear recurrences with polynomial coefficients satisfied by the sequences themselves, enjoy

{\tt https://sites.math.rutgers.edu/\~{}zeilberg/tokhniot/oGDW2.txt } \quad .

[Of course this file includes all the theorems in the previous file].

$\bullet$ 
If you want to see (one) theorem about the algebraic equation satisfied by the generating function enumerating generalized Dyck paths with set of steps $\{-4,-3,-2,-1,0,1,2,3,4\}$, look at

{\tt https://sites.math.rutgers.edu/\~{}zeilberg/tokhniot/oGDW3.txt } \quad . 

It has degree $41$ in $X(t)$!

{\bf  Enumeration According to the Sum of Areas of Generalized Dyck Paths}

$\bullet$ If you want to see $16$ theorems stating algebraic equations for the SUM OF THE AREAS under generalized Dyck paths for all non-trivial subsets of $\{-2,-1,0,1,2\}$,
as well, as estimates for the asymptotics of the average area (divided by $n^{3/2}$ (the ratio always tends to some constant)), see

{\tt https://sites.math.rutgers.edu/\~{}zeilberg/tokhniot/oGDW7.txt } \quad . 

$\bullet$ If you want to see interesting information about generalized Dyck paths for all non-trivial subsets of $\{-2,-1,0,1,2\}$,
that consist of linear recurrences (obtained by guessing, but definitely correct, since we know they exist from theoretical reasons), for
the straight enumeration, enumeration by `sum of areas', and sometimes, enumeration by `sum of area-squared', enabling estimates not only of
the asymptotic average area, but also of the asymptotic variance, look here:

{\tt https://sites.math.rutgers.edu/\~{}zeilberg/tokhniot/oGDW8.txt } \quad . 

{\bf Conclusion} 

We described an interesting application of Bruno Buchberger's seminal Gr\"obner Basis algorithm, and at the same time
demonstrated how an important class of enumerative combinatorics problems can be {\bf fully automated}, using the great
power of modern Computer Algebra Systems (in our case {\tt Maple}),  also using  calculus (that Maple knows very well), 
as well as good-old numerical (and symbolic)
dynamical programming. We also demonstrated how to teach the computer to build, {\it ab initio}, the system of equations (sometimes quite large)
before asking it  to kindly solve it.

{\bf References}

[AlZ] Gert Almkvist and Doron Zeilberger, {\it The method of differentiating under the integral sign},
J. Symbolic Computation {\bf 10} (1990), 571-591. \hfill\break
{\tt https://sites.math.rutgers.edu/\~{}zeilberg/mamarim/mamarimPDF/duis.pdf} \quad. \hfill\break
[It is implemented in procedure {\tt AZd} in the Maple package \hfill\break
{\tt https://sites.math.rutgers.edu/\~{}zeilberg/tokhniot/EKHAD.txt} ]

[AyZ] Arvind Ayyer and Doron Zeilberger, {\it The Number of [Old-Time] Basketball games with Final Score n:n where the Home Team was never losing but also never ahead by more than w Points},
Electronic J. of Combinatorics {\bf 14(1)} (2007), R19 [8 pp]. \hfill\break
{\tt https://sites.math.rutgers.edu/\~{}zeilberg/mamarim/mamarimhtml/basketball.html} \quad .

[BKKKKNW] C. Banderier, C. Krattenthaler, A. Krinik, D. Kruchinin, V. Kruchinin, D. Nguyen, and M. Wallner, 
{\it Explicit formulas for enumeration of lattice paths: basketball and the kernel method}, arXiv preprint arXiv:1609.06473 [math.CO], 2016. \hfill\break
https://arxiv.org/abs/1609.06473

[B1] Bruno Buchberger, {\it ``Ein Algorithmus zum Auffinden der Basiselemente des Restklassenringes nach einem nulldimensionalen Polynom ideal''},
PhD thesis, University of Innsbruck, 1965.

[B2] Bruno Buchberger, {\it Gr\"obner Bases: an algorithmic  method in polynomial idea theory}. In:
{\it ``Multidimensional  Systems Theory''}, ed. N.K. Bose, {\it ``Mathematics and Its Applications''}, chapter {\bf 6}, 184-232. D. Reidel Publishing company,
Dordrecht, 1985.

[Ch] Robin Chapman, {\it Moments of Dyck paths}, Discrete Mathematics {\bf 204} (1999), 113-117.

[CLO] David Cox, John Little, and Donald O'Shea, {\it ``Ideals, Varieties, and Algorithms''}, Springer, 1991.

[D] Robert Dougherty-Bliss, {\it Integral Recurrences from A to Z},   American Mathematical Monthly {\bf 129} (2022), 805-815. \hfill\break
{\tt https://arxiv.org/abs/2102.10170}

[Ek1] Bryan Ek, {\it ``Unimodal Polynomials and Lattice Walk Enumeration with Experimental Mathematics''}, 
PhD thesis, Rutgers University, May 2018. Available from \hfill\break
{\tt http://sites.math.rutgers.edu/\~{}zeilberg/Theses/BryanEkThesis.pdf} \quad .

[Ek2] Bryan Ek, {\it Lattice Walk Enumeration}, 29 March, 2018, {\tt https://arxiv.org/abs/1803.10920}.

[EkhZ] Shalosh B. Ekhad and Doron Zeilberger, {\it The Method(!) of ``Guess and Check''}, The Personal journal of
Shalosh B. Ekhad and Doron Zeilberger, Feb. 15, 2015. \hfill\break
{\tt https://sites.math.rutgers.edu/\~{}zeilberg/mamarim/mamarimhtml/gac.html} \quad . \hfill\break

[KP] Manuel Kauers and Peter Paule, {\it ``The Concrete Tetrahedron}, Springer, 2011.

[MSV] R. Merlini, R. Sprugnoli, and M.C. Verri, {\it The area determined by under-diagonal paths}, Proc. CAAP, vol. {\bf 1059},
Lecture Notes in Computer Science, 1996, 59-71.

[Sl] Neil A. J. Sloane, {\it The On-Line Encyclopedia of Integer Sequences® (OEIS)}, {\tt https://oeis.org/}.

[SRW] L.W. Shapiro, D. Rogers, and W. -J. Woan, {\it The Catalan numbers, the Lebesgue integral, and $4^{n-2}$}, Amer. Math. Monthly {\bf 104} (1997). 926-931.

[TZ] Thotsaporn Aek Thanatipanonda and Doron Zeilberger, {\it A Multi-Computational Exploration of Some Games of Pure Chance},
J. of Symbolic Computation {\bf 104}(2021), 38-68. \hfill\break
{\tt https://sites.math.rutgers.edu/\~{}zeilberg/mamarim/mamarimhtml/chance.html}

[Wo] W.-J. Woan, {\it Area of Catalan Paths}, Discrete Math. {\bf 226} (2001), 439-444.

\bigskip
\hrule
\bigskip
AJ Bu and Doron Zeilberger, Department of Mathematics, Rutgers University (New Brunswick), Hill Center-Busch Campus, 110 Frelinghuysen
Rd., Piscataway, NJ 08854-8019, USA. \hfill\break
Email: {\tt   ab1854  at math dot rutgers dot edu} \quad, \quad {\tt DoronZeil at gmail dot com}   \quad .

Written: {\bf May  15, 2023}.

\end